\documentclass{amsart}

\usepackage[frame,cmtip,arrow,matrix,line,graph,curve]{xy}
\usepackage{amsmath,amssymb}
\usepackage{latexsym}

\numberwithin{equation}{section}

\newtheorem{theorem}{Theorem}[section]
\newtheorem{proposition}[theorem]{Proposition}
\newtheorem{corollary}[theorem]{Corollary}
\newtheorem{lemma}[theorem]{Lemma}
\newtheorem{conjecture}[theorem]{Conjecture}
\newtheorem{remit}[theorem]{Remark}

\newenvironment{remark}{\begin{remit}\rm}{\end{remit}}

\newcommand{\pp}{\mathbb{P}}

\newcommand{\cc}{\mathbb{C}}
\newcommand{\rr}{\mathbb{R}}
\newcommand{\zz}{\mathbb{Z}}

\newcommand{\Hom}{\mathrm{Hom}}

\newcommand{\git}{/\!\!/}

\newcommand{\cA}{\mathcal{A} }
\newcommand{\cB}{\mathcal{B} }
\newcommand{\cF}{\mathcal{F} }
\newcommand{\cE}{\mathcal{E} }
\newcommand{\cN}{\mathcal{N} }
\newcommand{\cL}{\mathcal{L} }
\newcommand{\cO}{\mathcal{O} }
\newcommand{\cP}{\mathcal{P} }
\newcommand{\End}{\mathrm{End} }
\newcommand{\tU}{\tilde{U} }
\newcommand{\tcN}{\tilde{\cN} }
\newcommand{\tcE}{\tilde{\cE} }
\newcommand{\tcF}{\tilde{\cF} }
\newcommand{\cD}{\mathcal{D} }
\newcommand{\tV}{\tilde{V} }
\newcommand{\tcD}{\tilde{\cD} }
\newcommand{\tD}{\tilde{D} }
\newcommand{\tY}{\tilde{Y} }

\newcommand{\tdel}{\tilde{\Delta} }

\newcommand{\e}{\epsilon }
\newcommand{\s}{\sigma }
\newcommand{\g}{\gamma }
\newcommand{\hs}{\hat{\sigma} }
\newcommand{\he}{\hat{\epsilon} }
\newcommand{\hg}{\hat{\gamma} }

\begin{document}

\title[Desingularizations of moduli space]{Desingularizations of
the moduli space of rank 2 bundles over a curve}
\date{}
\author{Young-Hoon Kiem}
\address{Dept of Mathematics, Seoul National University,
Seoul, 151-747, Korea} \email{kiem@math.snu.ac.kr}
\author{Jun Li}
\address{Department of Mathematics, Stanford University, Stanford,
USA}\email{jli@math.stanford.edu}
\thanks{Young-Hoon Kiem was partially supported by KOSEF and SNU;
Jun Li was partially supported by NSF grants.} 


\begin{abstract} Let $X$ be a smooth projective curve of genus
$g\ge 3$ and $M_0$ be the moduli space of rank 2 semistable
bundles over $X$ with trivial determinant. There are three
desingularizations of this singular moduli space constructed by
Narasimhan-Ramanan \cite{NR}, Seshadri \cite{Se1} and Kirwan
\cite{k5} respectively. The relationship between them has not been
understood so far. The purpose of this paper is to show that there
is a morphism from Kirwan's desingularization to Seshadri's, which
turns out to be the composition of two blow-downs. In doing so, we
will show that the singularities of $M_0$ are terminal and the
plurigenera are all trivial. As an application, we compute the
Betti numbers of the cohomology of Seshadri's desingularization in
all degrees. This generalizes the result of \cite{BaSe} which
computes the Betti numbers in low degrees. Another application is
the computation of the stringy E-function (see \cite{Bat} for
definition) of $M_0$ for any genus $g\ge 3$ which generalizes the
result of \cite{kiem}.
\end{abstract}
 \maketitle

\begin{center} \emph{Dedicated to Professor Ronnie Lee.}\end{center}



\section{Introduction}

Let $X$ be a smooth projective curve of genus $g\ge 3$. Let $M_0$
be the moduli space of rank 2 semistable bundles over $X$ with
trivial determinant, which is a singular projective variety of
dimension $3g-3$. There are three desingularizations of $M_0$.
\begin{enumerate}
\item Seshadri's desingularization $S$ : fine moduli space of
parabolic bundles of rank 4 and degree zero such that the
endomorphism algebra of the underlying vector bundle is isomorphic
to a specialization of the matrix algebra $M(2)$. This is
constructed in \cite{Se1}.
\item Narasimhan-Ramanan's desingularization $N$ : moduli space of
Hecke cycles, as an irreducible subvariety of the Hilbert scheme
of conics. This is constructed in \cite{NR}.
\item Kirwan's desingularization $K$ : the result of systematic
blow-ups of $M_0$, constructed in \cite{k5}.
\end{enumerate}
For cohomological computation, $K$ is most useful thanks to the
Kirwan theory \cite{k2,KirL,k5}. On the other hand, $S$ and $N$
are moduli spaces themselves. The relationship between these
desingularizations has not been understood.

The first main result of this paper is that there is a birational
morphism (Theorem \ref{th4.1})
$$\rho:K\to S .$$
Since both $S$ and $K$ contain the open subset $M^s_0$ of stable
bundles, there is a rational map $\rho':K\dashrightarrow S$. By
GAGA and Riemann's extension theorem \cite{Kran}, it suffices to
show that $\rho'$ can be extended to a continuous map with respect
to the usual complex topology. By Luna's slice theorem, for each
point $x\in M_0-M_0^s$, there is an analytic submanifold $W$ of
the Quot scheme whose quotient by the stabilizer $H$ of a point in
both $W$ and the closed orbit represented by $x$ is analytically
equivalent to a neighborhood of $x$ in $M_0$. Furthermore,
Kirwan's desingularization $\tilde{W}\git H$ of $W\git H$ is a
neighborhood of the preimage of $x$ in $K$ by construction. Our
strategy is to construct a nice family of (parabolic) vector
bundles of rank 4 parametrized by $\tilde{W}$, starting from the
family of rank 2 bundles parametrized by $W$, which is induced
from the universal bundle over the Quot scheme. This is achieved
by successive applications of elementary modifications. Because
$S$ is the fine moduli space of such parabolic bundles of rank 4,
we get a morphism $\tilde{W}\to S$. This is invariant under the
action of $H$ and hence we have a morphism $\tilde{W}\git H\to S$.
Therefore, $\rho'$ extends to a neighborhood of the preimage of
$x$ in $K$.

The second main result of this paper is that the above morphism
$\rho$ is in fact the consequence of two blow-downs which can be
described quite explicitly (Theorem \ref{th5.6}). To prove this
theorem, we first show that Kirwan's desingularization $K$ can be
blown down twice by finding extremal rays. O'Grady in
\cite{ogrady} worked out such contractions for the moduli space of
rank 2 sheaves on a K3 surface. Since the proofs are almost same
as his case, we provide only the outline and necessary
modifications in \S \ref{contractions}. Next, we show that $\rho$
is constant along the fibers of the blow-downs and thus $\rho$
factors through the blown-down of $K$. Finally, Zariski's main
theorem tells us that $S$ is isomorphic to the blown-down. Using
this theorem, we can compute the discrepancy divisor of
$\pi_K:K\to M_0$ (Proposition \ref{discr}) and show that the
singularities are terminal. This implies that the plurigenera of
$M_0$ (or $K$, or $S$) are all trivial (Corollary \ref{cor5.4}).
We conjecture that the intermediate variety between $K$ and $S$ is
the desingularization $N$ by Narasimhan and Ramanan.

Our third main result is the computation of the cohomology of $S$.
In \cite{Ba,BaSe}, Balaji and Seshadri provides an algorithm for
the Betti numbers of $S$ for degrees up to $2g-4$. The cohomology
of Kirwan's \emph{partial} desingularization is computed in
\cite{k5} and $K$ is obtained as a single blow-up of this partial
desingularization. Since it is well-known how to compare
cohomology groups after blow-up (or blow-down) along a smooth
submanifold of an orbifold (\cite{GH} p.605), we can compute the
cohomology of $S$.

The last result of this paper is the computation of the stringy
E-function of $M_0$. The stringy E-function is a new invariant of
singular varieties, obtained as the measure of the arc space (see,
for instance, \cite{Bat}). From the knowledge of the discrepancy
divisor (Proposition \ref{discr}) and explicit descriptions of the
exceptional divisors of $\pi_K:K\to M_0$ (Proposition
\ref{prop5.1}), we show that
$$\begin{array}{ll}E_{st}(M_0)=&
\frac{(1-u^2v)^g(1-uv^2)^g-(uv)^{g+1}(1-u)^g(1-v)^g}{(1-uv)(1-(uv)^2)}\\
&- \frac{(uv)^{g-1}}{2}\big(
\frac{(1-u)^g(1-v)^g}{1-uv}-\frac{(1+u)^g(1+v)^g}{1+uv} \big).
\end{array}$$
Surprisingly, this is equal to the E-polynomial of the
intersection cohomology of $M_0$ when $g$ is even. For $g$ odd,
$E_{st}(M_0)$ is not a polynomial. As a consequence, the stringy
Euler number is
$$e_{st}(M_0):=\lim_{u,v\to
1}E_{st}(M_0)=4^{g-1}.$$
 If we denote by $e_g$ the stringy Euler
number of the moduli space $M_0$ for a genus $g$ curve, then the
equality
$$\sum_{g}e_gq^g=\frac14\frac1{1-4q}$$
holds for degree $\ge 2$. The coefficient $\frac14$ might be
related to the ``mysterious" coefficient $\frac14$ for the
S-duality conjecture test in \cite{VW}.

This paper is organized as follows. In sections 2 and 3, we review
Seshadri's and Kirwan's desingularizations respectively. In
section 4, we construct a morphism $\rho:K\to S$ by elementary
modification. In section 5, we show that $\rho$ is the composition
of two blow-downs. In section 6, we compute the cohomology of $S$
and the stringy E-function of $M_0$.

The first named author thanks Professor Ramanan for useful
conversations at the Korea Institute for Advanced Study during the
spring of 2003 and Professor Ronnie Lee for illuminating
discussions at Yale on desingularizations about 5 years ago. Part
of this paper was written while the first named author was
visiting Stanford University and Fudan University. Their
hospitality is greatly appreciated.


\section{Seshadri's desingularization}

Let $X$ be a compact Riemann surface of genus $g\ge 3$. Let
$M_0=M_X(2,\cO)$ denote the moduli space of semistable vector
bundles over $X$ of rank 2 with trivial determinant. Then $M_0$ is
a \emph{singular} normal projective variety of (complex) dimension
$3g-3$. In \cite{Se1}, Seshadri constructed a desingularization
$$\pi_S:S\to M_0$$
which restricts to an isomorphism on $\rho_S^{-1}(M_0^s)$ where
$M_0^s$ denotes the open subset of stable bundles. In fact, this
is constructed as the fine moduli space of a moduli problem which
we recall in this section. The main reference is \cite{Se2}
Chapter 5 and \cite{BaSe}.

Fix a point $x_0\in X$. Let $E$ be a vector bundle of rank 4 and
degree 0 on $X$ and $0\ne s\in E^*_{x_0}$ be a parabolic structure
with parabolic weights $0<a_1<a_2<1$.

\begin{lemma}\label{lem1} (\cite{Se2} 5.III Lemma 5) There are real
numbers $a_1, a_2$ such that for any semistable parabolic bundle
$(E,s)$ of rank 4 and degree 0, we have
\begin{enumerate}
\item $(E,s)$ is stable
\item $E$ is a semistable vector bundle.
\end{enumerate}
\end{lemma}

If we take sufficiently small $a_1$ and $a_2$, it is easy to see
that the conditions of the lemma are satisfied. Let us fix such
$a_1, a_2$.

It is well-known from \cite{MS} that the moduli functor
\begin{equation}\label{eq1}
\mathcal{P}:\mathcal{V}ar\to \mathcal{S}ets \end{equation}
 which
assigns to each variety $T$ the set of equivalence classes of
families of stable parabolic bundles of rank 4 and degree 0 over
$X$ parameterized by $T$, is represented by a smooth projective
variety, which we denote by $P$. It turns out that Seshadri's
desingularization $S$ is a closed subvariety of $P$.

We need a few more facts from \cite{Se2}  (Chapter 5, Propositions
7, 8, 9).
\begin{proposition}\label{prop1}
Let $E$ be a semistable vector bundle of rank 4 and degree 0 on
$X$. There is $0\ne s\in E^*_{x_0}$ such that the parabolic bundle
$(E,s)$ is stable if and only if for  any line bundle $L$ on $X$
of degree 0 there is no injective homomorphism of vector bundles
$$L\oplus L\hookrightarrow E.$$
\end{proposition}

\begin{proposition}\label{prop2}
Let $(E,s)$ be a stable parabolic bundle of rank 4 and degree 0.
Then the algebra $\End E$ of endomorphisms of the underlying
vector bundle $E$ has dimension $\le 4$. Moreover, if the algebra
$\End E$ is isomorphic to the matrix algebra $M(2)$ of $2\times 2$
matrices, then $E\cong F\oplus F$ for a unique stable vector
bundle $F$ of rank 2 and degree 0.
\end{proposition}

\begin{proposition}\label{prop3}
Let $(E_1,s_1)$, $(E_2,s_2)$ be two stable parabolic bundles of
rank 4, degree 0 over $X$. Suppose $\dim \End E_1=\dim \End
E_2=4$. Then they are isomorphic as parabolic bundles if and only
if the underlying vector bundles $E_1$ and $E_2$ are isomorphic.
\end{proposition}

Let $S'$ be the subset of $P$ consisting of stable parabolic
bundles $(E,s)$ such that $\End E\cong M(2)$ and $\mathrm{det} E$
is trivial. Then Proposition \ref{prop2} says we have a map $S'\to
M_0^s$ from $S'$ to the set of stable vector bundles. By
Proposition \ref{prop3}, this map is injective. By Proposition
\ref{prop1}, it is surjective as well. Seshadri's
desingularization $S$ of $M_0$ is defined as the closure of $S'$
in $P$ which is nonsingular by \cite{BaSe} Proposition 1.
Furthermore, the morphism $S'\to M_0^s$ extends to a morphism
$\pi_S: S\to M_0$ such that for each $(E,s)\in S$, $\mathrm{gr} E
\cong F\oplus F$ where $F$ is the polystable bundle representing
the image of $(E,s)$ in $M_0$.

Fix a nonzero element $e_0\in
\cc^4$. Let $\cA(2)$ be the set of elements in $$\Hom
(\cc^4\otimes\cc^4, \cc^4)$$ which gives us an algebra structure
on $\cc^4$ with the identity element $e_0$. There is a subset of
$\cA(2)$ which consists of algebra structures on $\cc^4$,
isomorphic to the matrix algebra $M(2)$. Let $\cA_2$ be the
closure of this subset. An element of $\cA_2$ is called a
\emph{specialization} of $M(2)$. Obviously, there is a locally
free sheaf $W$ of $\cO_{\cA_2}$-algebras on $\cA_2$ such that for
every $z\in\cA_2$, $W_z\otimes \cc$ is the specialization of
$M(2)$ represented by $z$.

Let $\cF$ be the subfunctor of the functor $\cP$ \eqref{eq1}
defined as follows. For each variety $T$, $\cF(T)$ is the set of
equivalence classes of families $\cE\to T\times X$ of stable
parabolic bundles on $X$ of rank 4 and degree 0 that satisfies the
following property (*): \\
for any $t\in T$ there is a neighborhood $T_1$ of $t$ in $T$ and a
morphism $f:T_1\to \cA_2$ such that $f^*W\cong
(p_T)_*(\mathcal{E}nd \cE)|_{T_1}$ as $\cO_{T_1}$-algebras where
$p_T:T\times X\to T$ is the projection to $T$.

\begin{theorem} \label{th1} (\cite{Se2} Chapter 5, Theorem 15)
 The functor $\cF$ is represented by $S$.
\end{theorem}

The condition (*) can be weakened slightly by the following
proposition.
\begin{proposition}\label{prop4} (\cite{Se2} Chapter 5, Proposition 1)
Let $T$ be a complex manifold and $B$ be a holomorphic vector
bundle of rank 4 equipped with an $\cO_T$ algebra structure.
Suppose there is an open dense subset $T'$ of $T$ such that for
each $t\in T'$, $B_t\otimes \cc$ is a specialization of $M(2)$.
Then for every $t\in T$, there is a neighborhood $T_1$ of $t$ and
a morphism $f:T_1\to \cA_2$ such that $f^*W\cong B|_{T_1}$.
\end{proposition}

To prove this, it suffices to consider any open set of $T$ over
which $B$ is trivial. But in this trivial case, the proposition is
obvious.

The singular locus of $M_0$ is the Kummer variety $\mathfrak{K}$
or the complement of $M_0^s$, isomorphic to the quotient $Jac_0/
\zz_2$ of the Jacobian of degree 0 line bundles over $X$ by the
involution $L\to L^{-1}$. There are $2^{2g}$ fixed points
$\zz_2^{2g}=\{[L\oplus L^{-1}]\, :\, L\cong L^{-1}\}$ and we have
a stratification
\begin{equation}\label{decM} M_0=M_0^s\sqcup
(\mathfrak{K} -\zz_2^{2g})\sqcup \zz_2^{2g}.\end{equation}

On the other hand, Seshadri's desingularization $S$ is stratified
by the rank of the natural conic bundle on $S$ (\cite{Ba} \S3) and
thus we have a filtration by closed subvarieties
\begin{equation}\label{filtS}
S\supset S_1\supset S_2\supset S_3 \end{equation} such that
$S-S_1=\pi_S^{-1}(M_0^s)\cong M_0^s$.
\begin{proposition} (\cite{BaSe})
\begin{enumerate}
\item
The image $\pi_S(S_1-S_2)$ is precisely the middle stratum
$\mathfrak{K}-\zz_2^{2g}$. In fact, $S_1-S_2$ is a
$\pp^{g-2}\times \pp^{g-2}$ bundle over $\mathfrak{K}-\zz_2^{2g}$.
\item The image of $S_2$ is precisely
the deepest strata $\zz_2^{2g}$ and $S_2-S_3$ is the disjoint
union of $2^{2g}$ copies of a vector bundle of rank $g-2$ over the
Grassmannian $Gr(2,g)$ while $S_3$ is the disjoint union of
$2^{2g}$ copies of the Grassmannian $Gr(3,g)$.
\end{enumerate}\end{proposition}

We end this section with the following proposition which is the
key for our construction of the morphism from Kirwan's
desingularization to Seshadri's desingularization.
\begin{proposition}\label{key}
\begin{enumerate}
\item Let $\cE\to T\times X$ be a family of semistable holomorphic vector
bundles of rank 4 and degree 0 on $X$ parameterized by a complex
manifold $T$. Assume the following:
\begin{enumerate}
\item for any $t\in T$ and any line bundle $L$ of degree $0$ on $X$,
$L\oplus L$ is not isomorphic to a subbundle of $\cE|_{t\times X}$
\item there is an open dense subset $T'$ of $T$ such that $\End
(\cE|_{t\times X})\cong M(2)$ for any $t\in T'$.
\end{enumerate}
Then we have a holomorphic map $\tau: T\to S$.
\item Suppose a holomorphic map $\tau:T\to S$ is given.
Suppose $T$ is an open subset of a nonsingular quasi-projective
variety $W$ on which a reductive group $G$ acts such that every
point in $W$ is stable and the (smooth) geometric quotient $W/G$
exists. Furthermore, assume that there is an open dense subset
$W'$ of $W$ such that whenever $t_1,t_2\in T\cap W'$ are in the
same orbit, we have $\tau(t_1)=\tau(t_2)$. Then $\tau$ factors
through the (smooth) image $\overline{T}$ of $T$ in the quotient
$W/G$, i.e. we have a continuous map $\overline{T}\to S$ such that
the diagram
$$\xymatrix{
T \ar[rr]^{\tau}\ar[dr]&& S\\
&{\overline{T}}\ar[ur] }$$
%
commutes.
\end{enumerate}
\end{proposition}

\begin{proof}
(1) Let $E_t=\cE|_{t\times X}$. For each $t\in T$, there is a
parabolic structure $0\ne s_t\in (E_t)_{x_0}^*$ such that
$(E_t,s_t)$ is a stable parabolic bundle by (a) and Proposition
\ref{prop1}. Hence we get a set-theoretic map $\tau: T\to P$.
Moreover, by (b), a dense open subset of $T$ is mapped to $S'$ and
thus $\tau$ is actually a map into $S$. We show that this is in
fact holomorphic.

By Proposition \ref{prop2}, $\dim \End E_t\le 4$. Since $\dim \End
E_t$ is an upper semi-continuous function of $t$, $\{t\in T\,|\,
\dim \End E_t=4\}$ is a closed subset of $T$. But there is a dense
open subset in $T$ where $\dim \End E_t=4$ by (b). Hence, $\dim
\End E_t=4$ for all $t\in T$. Consequently,
$(p_T)_*\mathcal{E}nd(\cE)$ is a locally free sheaf of
$\cO_{T}$-algebras of rank 4.

Since stability is an open property, there is a neighborhood $T_1$
of $t$ and $s\in \cE|_{T_1\times x_0}$ such that $(E_{t'},
s_{t'})$ is a stable parabolic bundle for every $t'\in T_1$.
Therefore $(\cE|_{T_1\times X}, s)$ is a family of stable
parabolic bundles and $(p_{T_1})_*\mathcal{E}nd(\cE|_{T_1\times
X})$ is a locally free sheaf of $\cO_{T_1}$-algebras. Hence by
assumption (b) and Proposition \ref{prop4}, we see that
$(\cE|_{T_1\times X}, s)$ is a family of stable parabolic bundles
satisfying (*) above. By deformation theory, we have a linear map
from the tangent space of $T_1$ at $t'$ to the deformation space
of $(E_{t'},s_{t'})$ which is isomorphic to the tangent space of
$P$. This is the derivative of $\tau$ at $t'$.  So we see that
$\tau$ is a holomorphic map from $T_1$ to $S$. Because we can find
a covering of $T$ by such open sets $T_1$, we deduce that $\tau$
is holomorphic.

(2) This is an easy consequence of the \'etale slice theorem. In
particular, the image $\overline{T}$ is an open subset of $W/G$ in
the usual complex topology.
%
\end{proof}



\section{Kirwan's desingularization}


In this section we recall Kirwan's desingularization from
\cite{k5}. We refer to \cite{kiem} for a very explicit description
of this desingularization process for the genus 3 case.

Note that we have the decomposition \eqref{decM}. The idea is to
blow up $M_0$ along the deepest strata $\zz_2^{2g}$ and then along
the proper transform of the middle stratum $\mathfrak{K}$. Let
$M_1$ denote the result of the first blow-up and $M_2$ the second
blow-up. Kirwan's \emph{partial} desingularization is the
projective variety $M_2$ which we have to blow up one more time to
get the \emph{full} desingularization $K$.

The moduli space $M_0$ is constructed as the GIT quotient of a
smooth quasi-projective variety $\mathfrak{R}$, which is a subset
of the space of holomorphic maps from the Riemann surface to the
Grassmannian $Gr(2,p)$ of $2$-dimensional quotients of $\cc^p$
where $p$ is a large even number, by the action of $G=SL(p)$. Over
each point in the deepest strata $\zz_2^{2g}$ there is a unique
closed orbit in $\mathfrak{R}^{ss}$. By deformation theory, the
normal space of the orbit at a point $h$, which represents
$L\oplus L^{-1}$ where $L\cong L^{-1}$, is
\begin{equation}\label{Znormal}
H^1(End_0(L\oplus L^{-1}))\cong H^1(\mathcal{O})\otimes
sl(2)\end{equation} where the subscript $0$ denotes the trace-free
part. According to Luna's slice theorem, there is a neighborhood
of the point $[L\oplus L^{-1}]$ with $L\cong L^{-1}$, homeomorphic
to $H^1(\mathcal{O})\otimes sl(2)\git SL(2)$ since the stabilizer
of the point $h$ is $SL(2)$ (\cite{k5} (3.3)). More precisely,
there is an $SL(2)$-invariant locally closed subvariety $W$ in
$\mathfrak{R}^{ss}$ containing $h$ and an $SL(2)$-equivariant
morphism $W\to H^1(\cO)\otimes sl(2)$, \'etale at $h$, such that
we have a commutative diagram
\begin{equation}\label{etalediag}
\xymatrix{ G\times_{SL(2)} \left(H^1(\cO)\otimes sl(2)\right)
\ar[d] & G\times_{SL(2)} W
\ar[d]\ar[r]\ar[l] &\mathfrak{R}^{ss}\ar[d] \\
H^1(\cO)\otimes sl(2)\git SL(2) & W\git SL(2)\ar[r]\ar[l] & M_0
}\end{equation} whose horizontal morphisms are all \'etale.

Next, we consider the middle stratum $\mathfrak{K}-\zz_2^{2g}$.
For each point, the normal space to the unique closed orbit over
it at a point $h$ representing $L\oplus L^{-1}$ with $L\neq
L^{-1}$, is isomorphic to
\begin{equation}\label{Knormal}
H^1(End_0(L\oplus L^{-1}))\cong H^1(\mathcal{O})\oplus
H^1(L^2)\oplus H^1(L^{-2}).\end{equation} The stabilizer $\cc^*$
acts with weights $0,2,-2$ respectively on the components. Hence,
there is a neighborhood of the point $[L\oplus L^{-1}]\in
\mathfrak{K}-\zz_2^{2g}$ in $M_0$, homeomorphic to
$$H^1(\mathcal{O})\bigoplus \big(H^1(L^2)\oplus
H^1(L^{-2})\git \cc^*  \big).$$ Notice that $H^1(\mathcal{O})$ is
the tangent space to $\mathfrak{K}$ and hence
$$
H^1(L^2)\oplus H^1(L^{-2}) \git \cc^* \cong \cc^{2g-2}\git\cc^*$$
is the normal cone. The GIT quotient of the projectivization $\pp
\cc^{2g-2}$ by the induced $\cc^*$ action is
$\pp^{g-2}\times\pp^{g-2}$ and the normal cone
$\cc^{2g-2}\git\cc^*$ is  obtained by collapsing the zero section
of the line bundle
$\mathcal{O}_{\pp^{g-2}\times\pp^{g-2}}(-1,-1)$.

Let $H$ be a reductive subgroup of $G=SL(p)$ and define $Z^{ss}_H$
as the set of semistable points in $\frak{R}^{ss}$ fixed by $H$.
Let $\frak{R}_1$ be the blow-up of $\frak{R}^{ss}$ along the
smooth subvariety $GZ^{ss}_{SL(2)}$. Then by Lemma 3.11 in
\cite{k2}, the GIT quotient $\frak{R}^{ss}_1\git G$ is the first
blow-up $M_1$ of $M_0$ along $GZ^{ss}_{SL(2)}\git G\cong
\zz_2^{2g}$. The $\cc^*$-fixed point set in $\frak{R}_1^{ss}$ is
the proper transform $\tilde{Z}^{ss}_{\cc^*}$ of $Z_{\cc^*}^{ss}$
and the quotient of $G\tilde{Z}^{ss}_{\cc^*}$ by $G$ is the
blow-up $\tilde{\mathfrak{K}}$ of $\mathfrak{K}$ along
$\zz_2^{2g}$. If we denote by $\frak{R}_2$ the blow-up of
$\frak{R}^{ss}_1$ along the smooth subvariety
$G\tilde{Z}^{ss}_{\cc^*}=G\times
_{N^{\cc^*}}\tilde{Z}^{ss}_{\cc^*}$ where $N^{\cc^*}$ is the
normalizer of $\cc^*$, the GIT quotient $\frak{R}_2^{ss}\git G$ is
the second blow-up $M_2$ again by Lemma 3.11 in \cite{k2}. This is
Kirwan's \emph{partial} desingularization of $M_0$ (See \S 3
\cite{k5}).

The points with stabilizer greater than the center $\{\pm 1 \}$ in
$\frak{R}_2^{ss}$ is precisely the exceptional divisor of the
second blow-up and the proper transform $\tilde{\Delta}$ of the
subset $\Delta$ of the exceptional divisor of the first blow-up,
which corresponds, via Luna's slice theorem, to $$SL(2)\pp \{\left(\begin{matrix} 0 & b\\
c&0\end{matrix}\right)\,|\, b,c\in H^1(\cO)\}\subset \pp
(H^1(\cO)\otimes sl(2)).$$ This is a simple exercise. Hence, if we
blow up $M_2$ along $\tilde{\Delta}\git SL(2)$, we get a smooth
variety $K$,  Kirwan's desingularization.


\section{Construction of the morphism}\label{constr}

The goal of this section is to prove the following.
\begin{theorem}\label{th4.1} There is a birational morphism
$$\rho:K\to S$$
from Kirwan's desingularization $K$ to Seshadri's
desingularization $S$.\end{theorem}

 Since the
desingularization morphisms
$$\pi_K:K\to M_0,\ \ \ \ \pi_S:S\to M_0$$
are both isomorphisms over $M^s_0$, we have a rational map
$$\rho':K \dashrightarrow S .$$
By GAGA (\cite{Hart} Appendix B, Ex.6.6), it suffices to find a
holomorphic map $\rho:K\to S$ that extends $\rho'$. By Riemann's
extension theorem \cite{Kran}, it suffices to show that $\rho'$
can be extended to a continuous map with respect to the usual
complex topology.

\subsection{Points over the middle stratum}\label{midsec}
Let us first extend to points over the middle stratum of $M_0$.
Let $l=[L\oplus L^{-1}]\in \mathfrak{K}-\zz_2^{2g}\subset M_0$ and
let $W_l$ be the \'etale slice of the unique closed orbit in
$\mathfrak{R}^{ss}$ over $l$. By Luna's slice theorem we have a
commutative diagram
\begin{equation}\label{etalediag2}
\xymatrix{ G\times_{\cc^*} \cN_l \ar[d] & G\times_{\cc^*} W_l
\ar[d]\ar[r]\ar[l] &\mathfrak{R}^{ss}\ar[d] \\
\cN_l\git \cc^* & W_l\git \cc^*\ar[r]\ar[l] & M_0 }\end{equation}
whose horizontal morphisms are all \'etale where $G=SL(p)$ and
$$\cN_l=H^1(\End(L\oplus L^{-1})_0)=H^1(\cO)\oplus H^1(L^2)\oplus H^1(L^{-2}) .$$
The slice $W_l$ is a subvariety of $\mathfrak{R}^{ss}$ and the
universal bundle over $\mathfrak{R}^{ss}\times X$ gives us a
vector bundle over $W_l\times X$. Since $W_l\to \cN_l$ is \'etale,
this gives us a holomorphic family $\cF$ of semistable vector
bundles over $X$ parametrized by a neighborhood $U_l$ of $0$ in
$\cN_l$. The idea now is to modify $\cF\oplus \cF$ to make it
satisfy the assumptions of Proposition \ref{key}.

The restriction of $\cF$ to $(U_l\cap H^1(\cO))\times X$ is a
direct sum $$\cL\oplus \cL^{-1}$$ where $\cL$ is a line bundle
coming from an \'etale map between $H^1(\cO)$ and the slice in the
Quot scheme for degree 0 line bundles.

To get Kirwan's desingularization, we blow up $\cN_l$ along
$H^1(\cO)$. Let $\pi_l:\tilde{\cN}_l\to \cN_l$ be the blow-up map.
Let $\tilde{U}_l=\pi_l^{-1}(U_l)\cap \tilde{\cN}^{ss}_l$ and $D_l$
be the exceptional locus in $\tilde{U}_l$. Let $\tilde{\cF}$ and
$\tilde{\cL}$ denote the pull-backs of $\cF$ and $\cL$ to
$\tilde{U}_l$ and $D_l$ respectively. Then we have surjective
morphisms
$$\tilde{\cF}|_{D_l}\to \tilde{\cL},\ \ \ \
\tilde{\cF}|_{D_l}\to \tilde{\cL}^{-1}.$$ Let $\tilde{\cF}'$  and
$\tilde{\cF}''$ be the kernels of
$$\tilde{\cF}\to \tilde{\cF}|_{D_l}\to \tilde{\cL},\ \ \ \
\tilde{\cF}\to\tilde{\cF}|_{D_l}\to \tilde{\cL}^{-1}$$
respectively. Define $\cE=\tilde{\cF}'\oplus \tilde{\cF}''$ over
$\tilde{U}_l\times X$.

\begin{lemma}\label{midex}
The bundle $\cE$ is a family of semistable vector bundles of rank
4 and degree 0 over $X$ parametrized by $\tilde{U}_l$ such that
the assumptions of Proposition \ref{key} are satisfied, i.e.
\begin{enumerate}
\item For each $t\in \tilde{U}_l$ and $L'\in Pic^0(X)$,
$L'\oplus L'$ is not isomorphic to any subbundle of $\cE|_{t\times
X}$.
\item $\cE|_{(\tilde{U}_l-D_l)\times X}\cong (\tilde{\cF}\oplus \tilde{\cF})|_{(\tilde{U}_l-D_l)\times X}$
and there is an open dense subset of $\tilde{U}_l$ where
$\mathrm{End}(\cE|_{t\times X})$ is a specialization of $M(2)$.
\item With respect to the action of $\cc^*$ on $\tilde{\cN}_l-D_l$, if
$t_1,t_2\in \tilde{U}_l-D_l$ are in the same orbit, then
$\cE|_{t_1\times X}\cong \cE|_{t_2\times X}$.
\end{enumerate}
\end{lemma}

\begin{proof}
Since $D_l$ is a smooth divisor in $\tilde{U}_l$, $\cE$ is locally
free of rank 4. Let $(a,b,c)\in \cN_l=H^1(\cO)\oplus
H^1(L^2)\oplus H^1(L^{-2})$. The weights of the $\cc^*$ action are
$0,2,-2$ respectively. It is well-known (see \cite[(2.5)
(iv)]{k5}) that the bundle $\cF|_{(a,b,c)\times X}$ is stable if
and only if the image of $(a,b,c)$ in $\mathfrak{R}^{ss}$ is a
stable point. This is equivalent to saying that $(a,b,c)$ is
stable with respect to the $\cc^*$ action. Hence
$\cF|_{(a,b,c)\times X}$ is stable if and only if $b\ne 0$ and
$c\ne 0$.

Let $t_0\in \tU_l-D_l$ and $\pi_l(t_0)=(a,b,c)$. This point has
nothing to do with the blow-up and the Hecke modification. Hence
$\tcE|_{t_0\times X}\cong \cF\oplus \cF|_{\pi_l(t_0)\times X}.$
The unstable points in $\tcN_l$ are the proper transform of
$\{(a,b,c)|b=0 \text{  or } c=0\}.$ Since $t_0$ is (semi)stable,
we have $b\ne 0$ and $c\ne 0$ which implies that
$F=\cF|_{\pi_l(t_0)\times X}$ is stable. Therefore, $\End(F\oplus
F)\cong M(2)$ which proves (2).

For $t_1, t_2\in \tU_l-D_l$, $\tcE|_{t_j\times X}\cong \cF\oplus
\cF|_{\pi_l(t_j)\times X}$ ($j=1,2$). But $\cF|_{\pi_l(t_1)\times
X}\cong \cF|_{\pi_l(t_2)\times X}$ if and only if $\pi_l(t_1)$ and
$\pi_l(t_2)$ are in the same orbit. This is equivalent to $t_1$
and $t_2$ being in the same orbit since $\tU_l-D_l$ is isomorphic
to the stable part of $\cN_l$. So we proved (3).

Let us prove (1). For $t\in \tU_l-D_l$, it is trivial since
$\tilde{\cF}'|_{t\times X}\cong \tilde{\cF}|_{t\times X}\cong
{\cF}|_{\pi_l(t)\times X}$ which is stable and the same is true
for $\tilde{\cF}''$.

 Let $C$ be a line in $\cN_l$ given by a map $\cc\to \cN_l$ with
$z\to (a, zb,zc)$ for $a\in H^1(\cO), 0\ne b\in H^1(L^2), 0\ne
c\in H^1(L^{-2})$. Note that any point in $D_l$ is represented by
such a line. Let $t$ be the point in $D_l$ represented by $C$.

 Let
$C_0=C\cap U_l$. By restricting $U_l$ if necessary, we can find an
open covering $\{V_i\}$ of $X$ such that $\cF|_{C_0\times V_i}$
are all trivial. Fix a trivialization for each $i$ and let
$L_a=\cL|_{a\times X}$. Since $\cF|_{0\times X}\cong L_a\oplus
L_a^{-1}$, the transition matrices are of the form
$$\left(\begin{matrix}
\lambda_{ij} & zb_{ij}\\
zc_{ij} & \lambda^{-1}_{ij}\end{matrix}\right)$$ where
$\lambda_{ij}|_{z=0}$ is the transition for $L_a$. The cocycle
condition tells us that
$$\{\lambda_{ij}b_{ij}|_{z=0}\},\ \ \ \ \ \{\lambda_{ij}^{-1}c_{ij}|_{z=0}\}
$$ are cocycles whose cohomology classes are nonzero because
$\cF|_{(a,zb,zc)\times X}$ is stable for $z\ne 0$. Let $\cF'$ be
the kernel of $\cF|_{C_0\times X}\to \cF|_{0\times X}\cong
L_a\oplus L_a^{-1}\to L_a$ where the first morphism is the
restriction and the last is the projection. Define $\cF''$ as the
kernel of $\cF|_{C_0\times X}\to \cF|_{0\times X}\cong L_a\oplus
L_a^{-1}\to L^{-1}_a$. Let $F'=\cF'|_{0\times X}$ and
$F''=\cF''|_{0\times X}$. Then by construction, $\tcF'|_{t\times
X}\cong F'$ and $\tcF''|_{t\times X}\cong F''$.

Any section of $\cF'$ over $C_0\times V_i$ is of the form $(zs_1,
s_2)$. Because
$$\left(\begin{matrix} s_1\\ s_2\end{matrix}\right) \longleftrightarrow
\left(\begin{matrix} zs_1\\ s_2\end{matrix}\right)\longmapsto
\left(\begin{matrix}
\lambda_{ij} & zb_{ij}\\
zc_{ij} &\lambda_{ij}^{-1}\end{matrix} \right)
\left(\begin{matrix} zs_1\\
s_2\end{matrix}\right)=
\left(\begin{matrix} z(\lambda_{ij}s_1+b_{ij}s_2)\\
\lambda_{ij}^{-1}s_2+z^2c_{ij}s_1\end{matrix}\right)
\longleftrightarrow \left(\begin{matrix}
\lambda_{ij}s_1+b_{ij}s_2\\
\lambda_{ij}^{-1}s_2+z^2c_{ij}s_1\end{matrix}\right)$$ the
transition for $\cF'$ is
$$\left(\begin{matrix}
\lambda_{ij} & b_{ij}\\
z^2c_{ij} & \lambda_{ij}^{-1}\end{matrix}\right) .$$ Hence $F'$
fits into a short exact sequence
$$0\to L_a\to F'\to L_a^{-1}\to 0$$
whose extension class is given by $\{\lambda_{ij}b_{ij}|_{z=0}\}$
which is nonzero. Hence, $F'=\cF'|_{z=0}$ is a nonsplit extension
of $L_a^{-1}$ by $L_a$ and similarly $F''=\cF''|_{z=0}$ is a
nonsplit extension of $L_a$ by $L_a^{-1}$. It is now an elementary
exercise to show that $E=F'\oplus F''$ does not have a subbundle
isomorphic to $L'\oplus L'$ for any $L'\in Pic^0(X)$. So we proved
(1).
\end{proof}

By Proposition \ref{key}, we have a holomorphic map from the image
of $\tilde{U}_l$ in $\tcN_l^{ss}/\cc^*$ to $S$.  Since the image
is open  in the usual complex topology by the slice theorem, this
implies that $\rho'$ extends continuously to a neighborhood of the
points in $K$ lying over $l$. Since $\rho'$ is defined on an open
dense subset, there is at most one continuous extension.
Therefore, the extensions for various points $l$ in the middle
stratum $\mathfrak{K}-\zz_2^{2g}$ are compatible and so $\rho'$ is
extended to all the points in $K$ except those over the deepest
strata $\zz_2^{2g}$.

\subsection{Points over the deepest strata}\label{setup}

Let us next extend $\rho'$ to the points over the deepest strata
$\zz_2^{2g}$.  The exactly same argument applies to all the points
in $\zz_2^{2g}$, so we consider only the points in $K$ over
$0=[\cO\oplus \cO]$. Let $W$ be the \'etale slice of the unique
closed orbit in $\mathfrak{R}^{ss}$ over $[\cO\oplus \cO]\in M_0$.
Let
$$\cN=H^1(\cO)\otimes sl(2).$$  By Luna's slice theorem, a
neighborhood of $[\cO\oplus \cO]$ in $M_0$ is analytically
equivalent to a neighborhood of the vertex $\overline{0}$ in the
cone $\cN\git SL(2)$ from the diagram \eqref{etalediag}. Hence a
neighborhood of the preimage of $[\cO\oplus \cO]$ in $K$ is
biholomorphic to an open set of the desingularization $\tcN\git
SL(2)$, obtained as a result of three blow-ups from $\cN\git
SL(2)$, described below. Therefore it suffices to construct a
holomorphic map from a neighborhood $\tV$ of the preimage of
$\overline{0}$ in $\tcN\git SL(2)$ to $S$.

Let $\Sigma$ be the subset of $\cN$ defined by
$$SL(2)\{H^1(\cO)\otimes \left(\begin{matrix}
1&0\\0&-1\end{matrix}\right)\}.$$ Let $\pi_1:\cN_1\to \cN$ be the
first blow-up in the partial desingularization process, i.e. the
blow-up at 0, and let $\cD^{(1)}_1$ be the exceptional divisor.
Recall that $\Delta$ is the subset of $\cD^{(1)}_1$ defined as
$$SL(2)\pp \{\left(\begin{matrix} 0 & b\\
c&0\end{matrix}\right)\,|\, b,c\in H^1(\cO)\}.$$ Let
$\tilde{\Sigma}$ be the proper transform of $\Sigma$ in $\cN_1$.
Then the singular locus of $\cN^{ss}_1\git SL(2)$ is the quotient
of $\Delta\cup \tilde{\Sigma}$ by $SL(2)$. It is an elementary
exercise to check that
\begin{equation}\label{eq4.0}\cD^{(1)}_1\cap
\tilde{\Sigma}=SL(2)\pp \{H^1(\cO)\otimes \left(\begin{matrix}
1&0\\0&-1\end{matrix}\right)\}=\Delta\cap
\tilde{\Sigma}.\end{equation}

 Let $\pi_2:\cN_2\to
\cN_1$ be the second blow-up, i.e. the blow-up along
$\tilde{\Sigma}$ and let $\cD^{(2)}_2$ be the exceptional divisor.
Let $\cD^{(1)}_2$ be the proper transform of $\cD^{(1)}_1$. The
singular locus of $\cN_2\git SL(2)$ is the quotient of the proper
transform $\tilde{\Delta}$ of $\Delta$.

Finally let $\pi_3:\tcN=\cN_3\to \cN_2$ denote the blow-up of
$\cN_2$ along $\tilde{\Delta}$ and let $\tcD^{(3)}=\cD^{(3)}_3$ be
the exceptional divisor while $\tcD^{(1)}=\cD^{(1)}_3$,
$\tcD^{(2)}=\cD^{(2)}_3$ are the proper transforms of
$\cD^{(1)}_2$ and $\cD^{(2)}_2$ respectively. Let $\pi:\tcN\to
\cN$ be the composition of the three blow-ups. Also let
$D^{(j)}_i$ be the quotient of $\cD^{(j)}_i$ in $\cN_i\git SL(2)$
for $1\le i\le 3$ and $1\le j\le i$.

As in the middle stratum case, the pull-back of the universal
bundle over $\mathfrak{R}^{ss}\times X$ to $W\times X$ gives us a
holomorphic family $\cF$ of rank 2 semistable vector bundles over
$X$ parametrized by an open neighborhood $U$ of $0$ in $\cN$. Let
$V$ be the image of $U$ under the good quotient morphism $\cN\to
\cN\git SL(2)$. Then $V$ is an open neighborhood of
$\overline{0}$. Let $U_1=\pi_1^{-1}(U)\cap \cN_1^{ss}$ and $V_1$
be the image of $U_1$ by the good quotient morphism $\cN_1\to
\cN_1\git SL(2)$. From the commutative diagram
$$\xymatrix{
\cN_1^{ss}\ar[r]\ar[d]_{\pi_1} &\cN_1\git SL(2)\ar[d]^{\overline{\pi}_1}\\
\cN\ar[r] & \cN\git SL(2)}$$ we see that
$V_1=\overline{\pi}_1^{-1}(V)$.

Let $U_2=\pi_2^{-1}(U_1)\cap \cN_2^{ss}$ and $V_2$ be the image of
$U_2$ in the quotient of $\cN_2$. Then we have
$V_2=\overline{\pi}_2^{-1}(V_1)$ where $\overline{\pi}_2:
\cN_2\git SL(2)\to \cN_1\git SL(2)$. Similarly, let
$\tU=\pi_3^{-1}(U_2)\cap \tcN^{ss}$ and $\tV$ be the image of
$\tU$ in the quotient of $\tcN$. By construction, $\tV$ is smooth
with simple normal crossing divisors $\tD^{(1)}, \tD^{(2)},
\tD^{(3)}$ where $\tD^{(j)}=D^{(j)}_3$. To simplify our notation
we denote the intersection of $\tD^{(2)}$ with $\tV$ again by
$\tD^{(2)}$.

Since we already extended $\rho'$ to the points over the middle
stratum, we have a holomorphic map $\rho':\tV-(\tD^{(1)}\cup
\tD^{(3)})\to S$ and we have to extend it to $\rho:\tV\to S$.

\subsection{Points in $\tD^{(1)}-(\tD^{(2)}\cup \tD^{(3)})$} \label{deep1}
In this subsection, we extend $\rho'$ to points in $\tV$ that lies
over the quotient of $\cD_1^{(1)}-\Delta$ via
$\overline{\pi}_3\circ\overline{\pi}_2$. Notice that
$\cD_1^{(1)}-\Delta$ does not intersect with the blow-up centers
of the second and third blow-up and hence it remains unchanged.

Our strategy is again to modify the pull-back of $\cF\oplus \cF$
to $U_1-\Delta\cup \tilde{\Sigma}$ so that $\rho'$ extends to a
holomorphic map near the quotient of $\cD_1^{(1)}-\Delta$ by
Proposition \ref{key}.

Let $\cF_1$ be the pull-back of $\cF$ to $U_1\times X$  via
$\pi_1\times 1_X$. Then $\cF_1|_{\cD^{(1)}_1\times X}\cong
\cO\oplus \cO$ since $\cF|_{0\times X}$ is trivial. Let $\cF'_1$
be the kernel of
$$\cF_1\to \cF_1|_{\cD^{(1)}_1\times X}\cong \cO_{\cD^{(1)}_1\times X}\oplus
\cO_{\cD^{(1)}_1\times X}\to \cO_{\cD^{(1)}_1\times X}$$ where the
second arrow is the projection onto the first component. Let
$\cF''_1$ be defined similarly with the projection onto the second
component. By computing transition matrices as in the proof of
Lemma \ref{midex}, it is immediate that $\cF'_1|_{t_1\times X}$
and $\cF''_1|_{t_1\times X}$ are nonsplit extensions of $\cO$ by
$\cO$ if $t_1=\left[\begin{matrix} a&b\\
c&-a\end{matrix}\right]\in \pp \cN=\cD^{(1)}_1$ with $b\ne 0$ and
$c \ne 0$ in $H^1(\cO)$.

Suppose $t_1\in \cD^{(1)}_1-\Delta$. Then $a,b,c$ are linearly
independent because otherwise we can find $g\in SL(2)$ such that
$gt_1g^{-1}$ is of the form
\begin{equation}\label{conju}\left[\begin{matrix} 0&*\\
*&0\end{matrix}\right]\ \ \ \ \text{ or }\ \ \ \
\left[\begin{matrix} *&0\\ *&*\end{matrix}\right]. \end{equation}
The first case belongs to $\Delta$ while the second is unstable in
$\tcN$ and is deleted after all. In particular,
 $a, b, c$ are all nonzero and thus $\cF'_1|_{t_1\times X}$ and
$\cF''_1|_{t_1\times X}$ are nonsplit extensions of $\cO$ by $\cO$
whose extension classes are $b, c$ respectively.

The inclusion $\cF_1'\hookrightarrow \cF_1$ gives us a
homomorphism $\cF_1'|_{\cD^{(1)}_1\times X}\to
\cF_1|_{\cD^{(1)}_1\times X}\cong \cO\oplus\cO$ whose image is the
second factor $\cO$ and the kernel of this homomorphism is $\cO$.
Similarly, the trivial bundle $\cO_{\cD_1^{(1)}\times X}$ is a
subbundle of $\cF''_1|_{\cD_1^{(1)}\times X}$ and we have a
diagonal embedding of $\cO_{\cD_1^{(1)}\times X}$ into
$\cF'_1\oplus \cF''_1|_{\cD_1^{(1)}\times X}$. Let $\cE_1$ be the
kernel of
$$\cF'_1\oplus \cF''_1\to \cF'_1\oplus \cF''_1|_{\cD^{(1)}_1\times
X}\to \cF'_1\oplus  \cF''_1|_{\cD^{(1)}_1\times
X}/\cO_{\cD_1^{(1)}\times X}.$$ As in the proof of Lemma
\ref{midex}, introduce a local coordinate $z$ of a suitable curve
passing through $t_0$ and write the transition for $\cF'_1\oplus
\cF''_1$ as \begin{equation} \left(\begin{matrix}
\lambda_{ij} & b_{ij}& 0& 0\\
z^2c_{ij} & \lambda_{ij}^{-1}&0&0\\
0&0&\lambda_{ij}&z^2b_{ij}\\
0&0&c_{ij}&\lambda_{ij}^{-1}
\end{matrix}\right)\end{equation}
where $\lambda_{ij}=1+za_{ij}$. Note that, when restricted to
$z=0$, the cocycles $\{a_{ij}\}$, $\{b_{ij}\}$, $\{c_{ij}\}$
represent the classes $a,b,c\in H^1(\cO)$ respectively.

A local section of $\cE_1$ as a subsheaf of $\cF'_1\oplus \cF''_1$
is of the form $(s_1, zs_2, zs_3, s_1+zs_4)$. Because
\begin{equation}\label{eq4.2}\begin{array}{ll}
\left(\begin{matrix} s_1\\
s_2\\s_3\\s_4\end{matrix}\right)&\leftrightarrow
\left(\begin{matrix} s_1\\
zs_2\\zs_3\\s_1+zs_4\end{matrix}\right)\mapsto
\left(\begin{matrix}
\lambda_{ij} & b_{ij}& 0& 0\\
z^2c_{ij} & \lambda_{ij}^{-1}&0&0\\
0&0&\lambda_{ij}&z^2b_{ij}\\
0&0&c_{ij}&\lambda_{ij}^{-1}
\end{matrix}\right)
\left(\begin{matrix} s_1\\
zs_2\\zs_3\\s_1+zs_4\end{matrix}\right) \\
&=\left(
\begin{matrix} \lambda_{ij}s_1+zb_{ij}s_2\\
z^2c_{ij}s_1+z\lambda_{ij}^{-1}s_2\\
z^2b_{ij}s_1+ z\lambda_{ij}s_3+z^3b_{ij}s_4\\
zc_{ij}s_3+\lambda_{ij}^{-1}s_1+z\lambda_{ij}^{-1}s_4
\end{matrix}\right) \leftrightarrow
\left(\begin{matrix} \lambda_{ij}s_1+zb_{ij}s_2\\
zc_{ij}s_1+\lambda_{ij}^{-1} s_2\\
zb_{ij}s_1+\lambda_{ij}s_3+z^2b_{ij}s_4\\
\frac{\lambda_{ij}^{-1}-\lambda_{ij}}{z}s_1-b_{ij}s_2+c_{ij}s_3+\lambda_{ij}^{-1}s_4
\end{matrix}\right),\end{array}\end{equation}
the transition for $\cE_1$ is
\begin{equation}\label{eq4.3}
\left(\begin{matrix} \lambda_{ij}&zb_{ij}&0&0\\
zc_{ij}&\lambda_{ij}^{-1}&0&0\\
zb_{ij}&0&\lambda_{ij}&z^2b_{ij}\\
-2a_{ij}&-b_{ij}&c_{ij}&\lambda^{-1}_{ij}
\end{matrix}\right).\end{equation}
Put $z=0$ to see that the transition for $\cE|_{t_1\times X}$ is
\begin{equation}\label{eq4.5}\left(\begin{matrix}
1&0&0&0\\0&1&0&0\\0&0&1&0\\-2a_{ij}|_{z=0}&-b_{ij}|_{z=0}&c_{ij}|_{z=0}&1\end{matrix}\right).
\end{equation}
Hence we have a filtration by subbundles
\begin{equation}\label{eq4.4}\cE|_{t_1\times X}=E_4\supset E_3\supset E_2\supset
E_1\supset E_0=0\end{equation} such that $E_{i+1}/E_i\cong \cO_X$.
The extension $E_2$ of $\cO$ by $E_1\cong \cO$ is nontrivial since
$c\ne 0$. An extension of $\cO$ by $E_2$ is parameterized by
$Ext^1(\cO,E_2)$ which fits in the exact sequence
$$\xymatrix{
Hom(\cO,\cO)\ar[r]^{c} &Ext^1(\cO,\cO)\ar[r]&Ext^1(\cO,E_2)\to
Ext^1(\cO,\cO)}$$ and $E_3$ is the image of $b\in
Ext^1(\cO,\cO)\cong H^1(\cO)$ which is nonzero since $b,c$ are
linearly independent. Hence $E_3$ is a nonsplit extension.
Similarly $E_4$ is a nonsplit extension since $a,b,c$ are linearly
independent. Hence \eqref{eq4.4} is the result of three nonsplit
extensions. This certainly implies that the condition (a) of
Proposition \ref{key} is satisfied for points in $\tU$ over
$\cD^{(1)}_1-\Delta$. The other conditions of Proposition
\ref{key} (1), (2) are trivially satisfied and hence $\rho'$
extends to the points over the quotient of the points over
$\cD^{(1)}_1-\Delta$ as desired.

\subsection{Points in  $\tD^{(3)}-\tD^{(2)}$}

We use the notation of \S \ref{deep1}. Suppose now
$t_1=\left[\begin{matrix} a&b\\
c&-a\end{matrix}\right]\in \Delta-\tilde{\Sigma}$. Then $a,b,c$
span 2-dimensional subspace of $H^1(\cO)$. The bundle
$\cE_1|_{t_1\times X}$ in the previous subsection  has transition
matrices of the form \eqref{eq4.5}. The one dimensional space of
linear relations of $a,b,c$ gives rise to an embedding of $\cO$
into $\cE_1|_{t_1\times X}$. More generally, the family of linear
relations of $a,b,c$ gives us a line bundle over
$\Delta-\tilde{\Sigma}$. Let $\cL_1$ denote the pull-back of this
line bundle to $(\Delta-\tilde{\Sigma})\times X$. Then we have an
embedding of $\cL_1$ into $\cE_1|_{(\Delta-\tilde{\Sigma})\times
X}$.
%
%
%
Let $\cE_3$ (resp. $\cL_3$)  be the pull-back of $\cE_1$ (resp.
$\cL_1$) to $\tilde{U}=U_3$ (resp. $\tcD^{(3)}-\tcD^{(2)}$).

 Let $\tilde{\cE}$ be the kernel of
$$\cE_3|_{}\to \cE_3|_{(\tcD^{(3)}-\tcD^{(2)})\times X}
\to \cE_3|_{(\tcD^{(3)}-\tcD^{(2)})\times X}/\cL_3.$$ We claim
that $\tcE$ satisfies the conditions of Proposition \ref{key} and
hence $\rho'$ extends to the quotient of $\tcD^{(3)}-\tcD^{(2)}$.

For simplicity, let $t_1$ be $\left[\begin{matrix} 0&b\\
c&0\end{matrix}\right]\in \Delta-\tilde{\Sigma}$ with $b,c$
linearly independent. (The general case is obtained by
conjugation.) Let $t_3\in \tcD^{(3)}-\tcD^{(2)}$ be a (semi)stable
point lying over $t_1$. Now we make local computations as in
\eqref{eq4.2} and \eqref{eq4.3}.

A point $t_3\in \tcD^{(3)}$ represents a normal direction to
$\Delta$ at $t_1$. Choose a local parameter $z$ of the direction
such that $z=0$ represents $t_1$.

If $t_3$ represents a normal direction of $\Delta$ tangent to
$\tcD^{(1)}$, then from \eqref{eq4.5}, the transition of the
restriction of $\cE_3$ to the direction is of the form
\begin{equation}\label{eq4.8}\left(\begin{matrix}
1&0&0&0\\0&1&0&0\\0&0&1&0\\-2zd_{ij}&-b_{ij}&c_{ij}&1\end{matrix}\right)
\end{equation}
for some cocycle $\{d_{ij}\}$ which gives rise to a nonzero class
$d\in H^1(\cO)$ at $z=0$ such that $d,b,c$ are linearly
independent. In this case, the transition for $\tcE|_{t_3\times
X}$ is of the form
\begin{equation}\label{eq4.10}\left(\begin{matrix}
1&0&0&0\\0&1&0&0\\0&0&1&0\\
-2d_{ij}|_{z=0}&-b_{ij}|_{z=0}&c_{ij}|_{z=0}&1\end{matrix}\right)\end{equation}
by a local computation. Hence, the condition (1) of Proposition
\ref{key} is satisfied because the bundle is obtained by three
nonsplit extensions.

 Suppose $t_3$ represents the
direction normal to $\cD^{(1)}$. Then we can use the same curve we
used in \S \ref{deep1} and the transition of $\cE_3$ is given by
\eqref{eq4.3}. More generally, the transition of $\cE_3$
restricted to the direction of any $t_3$, not tangent to
$\cD^{(1)}$, is of the form
\begin{equation}\label{eq4.7}\left(\begin{matrix}
1+za_{ij}&zb_{ij}&0&0\\
zc_{ij}&1-za_{ij}&0&0\\
zb_{ij}&0&1+za_{ij}&0\\
-2zd_{ij}&-b_{ij}&c_{ij}&1-za_{ij}\end{matrix}\right)\end{equation}
mod $z^2$ for some cocycle $\{d_{ij}\}$.
 A local section of $\tcE$ is of the form $(s_1,
zs_2, zs_3,zs_4)$ and by computing as in \eqref{eq4.2} starting
with \eqref{eq4.7}, we deduce that the transition for
$\tcE|_{t_3\times X}$ is of the form
\begin{equation}\label{eq4.9}\left(\begin{matrix}
1&0&0&0\\c_{ij}|_{z=0}&1&0&0\\b_{ij}|_{z=0}&0&1&0\\
-2d_{ij}|_{z=0}&-b_{ij}|_{z=0}&c_{ij}|_{z=0}&1\end{matrix}\right).
\end{equation}
This implies that the bundle has a filtration by subbundles as in
\eqref{eq4.4} obtained by three nonsplit extensions. Hence
$\tcE|_{t_3\times X}$ satisfies the condition (1) of Proposition
\ref{key}.

Because the other conditions of Proposition \ref{key} are
trivially satisfied on the stable part of $U$, we deduce that the
holomorphic map $\rho'$ extends to the quotient of
$\tU-\tcD^{(2)}$. So far, we extended $\rho'$ to the complement of
the quotient of $\tcD^{(2)}\cap (\tcD^{(1)}\cup \tcD^{(3)})$ which
consists of points lying over $\Delta\cap \tilde{\Sigma}$.

\subsection{Points in $\tD^{(2)}\cap (\tD^{(1)}\cup \tD^{(3)})$}
In this subsection, we finally extend $\rho'$ to everywhere in $K$
and finish the proof of Theorem \ref{th4.1}. We use the notation
of \S \ref{setup}. By the slice theorem, we have a map $\tV\to K$,
biholomorphic onto a neighborhood of the preimage of $[\cO\oplus
\cO]$. So it suffices to construct a holomorphic map $\tV\to S$.

We have a commutative diagram
$$\xymatrix{
\tV\ar@{^(->}[r]\ar[d]_{\alpha} & K\ar[d]^{\beta}\\
V_1\ar@{^(->}[r] &M_1}$$ where the vertical maps are blow-ups. We
already constructed a holomorphic map
$$\nu: \tV-\alpha^{-1}(\Delta\cap \tilde{\Sigma}\git SL(2))\to S$$

Let $x$ be any point in $\Delta\cap \tilde{\Sigma}\git SL(2)$.
%
From \eqref{eq4.0}, $x$ is represented by the orbit of
$\left[\begin{matrix} a^0&0\\ 0&-a^0\end{matrix}\right]$ for some
$[a^0]\in H^1(X, \cO)$. The stabilizer of the point in $SL(2)$ is
$\cc^*$ and the normal space $Y$ to its orbit is isomorphic to
$\cc^g\oplus \cc^{2g-2}$ where $\cc^g$ is the tangent space of the
blow-up $\widetilde{H^1(\cO)}=\mathrm{bl}_0H^1(\cO)$ and
$\cc^{2g-2}\cong H^1(\cO)/\cc a^0 \oplus H^1(\cO)/\cc a^0$.

Obviously, a neighborhood $Y_1$ of $0$ in $Y$ is holomorphically
embedded into $U_1$, perpendicular to the $SL(2)$-orbit of the
point $[a^0]$ and the vector bundle $\cF_1|_{Y_1\times X}$ has
transition matrices of the form
\begin{equation}\label{10.29}\left(\begin{matrix}
1+z_1(a^0_{ij}+a_{ij}) & z_1 b_{ij}\\
z_1c_{ij} &
1-z_1(a^0_{ij}+a_{ij})\end{matrix}\right).\end{equation} Here
$a=\{a_{ij}\}, b=\{b_{ij}\}, c=\{c_{ij}\}$ are classes in
$H^1(\cO)$, not parallel to $a^0$ if nonzero and $z_1$ is the
coordinate for the normal direction of $\pp H^1(\cO)$ in
$\widetilde{H^1(\cO)}$.

By Luna's \'etale slice theorem, a neighborhood of the vertex of
the cone $Y\git \cc^*$ is analytically equivalent to a
neighborhood of $x$ in $V_1$ or $M_1$. Let $\tilde{Y}$ denote the
proper transform of $Y_1$ in $\tU$. Then the image of $\tilde{Y}$
in $\tilde{V}$ is biholomorphic to a neighborhood of
$\alpha^{-1}(x)$. Our goal is to construct a family of rank 4
bundles on $X$ parametrized by $\tilde{Y}$ satisfying the
conditions of Proposition \ref{key}. Then we can conclude that
$\nu$ extends to $\alpha^{-1}(x)$.

Recall that we have a rank 2 bundle $\cF_1$ over $U_1\times X$.
Let $\cF_{Y_1}=\cF_1|_{Y_1\times X}$. Let $\cD_{Y_1}^{(1)}$ be the
divisor in $Y_1$ given by $z_1=0$. Then from \eqref{10.29} we see
that
$$\cF_{Y_1}|_{\cD_{Y_1}^{(1)}\times X}\cong \cO\oplus \cO.$$
Let $\cF_{Y_1}'$ (resp. $\cF_{Y_1}''$) be the kernel of
$$\cF_{Y_1}\to \cF_{Y_1}|_{\cD^{(1)}_{Y_1}\times X}\cong
\cO\oplus\cO\to \cO$$ where the last arrow is the projection onto
the first (resp. second) component. From a local computation as in
\S \ref{midsec}, the transition matrices of $\cF'_{Y_1}$ and
$\cF''_{Y_1}$ are respectively
$$\left(\begin{matrix}
1+z_1(a^0_{ij}+a_{ij}) & b_{ij}\\
z_1^2c_{ij} & 1-z_1(a^0_{ij}+a_{ij})\end{matrix}\right), \ \ \ \
\left(\begin{matrix}
1+z_1(a^0_{ij}+a_{ij}) & z_1^2b_{ij}\\
c_{ij} & 1-z_1(a^0_{ij}+a_{ij})\end{matrix}\right)
$$
In particular, $\cF_{Y_1}'$ and $\cF_{Y_1}''$ restricted to
$$\tilde{\Sigma}_{Y_1}=Y_1\cap \{b=c=0\}=Y_1\cap (\cc^g\oplus 0)\subset
\cc^g\oplus \cc^{2g-2}=Y$$ are given by transition matrices
$$\left(\begin{matrix}
1+z_1(a^0_{ij}+a_{ij}) & 0\\
0 & 1-z_1(a^0_{ij}+a_{ij})\end{matrix}\right) $$ and thus
$$\cF'_{Y_1}|_{\tilde{\Sigma}_{Y_1}\times X}\cong \cL_{Y_1}\oplus
\cL_{Y_1}^{-1}$$ for some line bundle $\cL_{Y_1}$ over
$\tilde{\Sigma}_{Y_1}\times X$.

Let $Y_2$ be the proper transform of $Y_1$ in $U_2$ by the blow-up
(and subtraction of unstable points) map $U_2\to U_1$. In other
words, $Y_2$ is the blow-up of $Y_1$ along $\tilde{\Sigma}_{Y_1}$
with unstable points removed. Let $z_2$ be the coordinate of the
normal direction of the exceptional divisor $\cD^{(2)}_{Y_2}$ at a
point $[b,c]$ over $(z_1, a)$. Let $\cF_{2,0}'$, $\cF_{2,0}''$ be
the pull-back of $\cF'_{Y_1}$, $\cF_{Y_1}''$ to $Y_2\times X$
respectively. Let $\cL_{Y_2}$ denote the pull-back of $\cL_{Y_1}$
to $\cD^{(2)}_{Y_2}\times X$.

Let $\cF_{Y_2}'$ be the kernel of
$$\cF_{2,0}'\to \cF_{2,0}'|_{\cD^{(2)}_{Y_2}\times X}\cong
\cL_{Y_2}\oplus\cL_{Y_2}^{-1}\to \cL_{Y_2}$$ and $\cF_{Y_2}''$ be
the kernel of
$$\cF_{2,0}''\to \cF_{2,0}''|_{\cD^{(2)}_{Y_2}\times X}\cong
\cL_{Y_2}\oplus\cL_{Y_2}^{-1}\to \cL_{Y_2}^{-1}.$$ Let
$\cD^{(1)}_{Y_2}$ be the proper transform of $\cD^{(1)}_{Y_1}$. By
a local computation, it is easy to see that the trivial bundle
$\cO$ is a subbundle of both $\cF_{Y_2}'|_{\cD^{(1)}_{Y_2}\times
X}$ and $\cF_{Y_2}''|_{\cD^{(1)}_{Y_2}\times X}$ as in \S 4.3. Let
$\cE_{Y_2}$ be the kernel of
$$\cF_{Y_2}'\oplus \cF_{Y_2}''\to
\cF_{Y_2}'\oplus\cF_{Y_2}''|_{\cD_{Y_2}^{(1)}\times X}\to
\cF_{Y_2}'\oplus \cF_{Y_2}''|_{\cD_{Y_2}^{(1)}\times X}/\cO.$$

 The
inclusion $\cE_{Y_2}\hookrightarrow \cF_{Y_2}'\oplus \cF_{Y_2}''$
induces $\cE_{Y_2}|_{\cD^{(1)}_{Y_2}\times X}\to \cF_{Y_2}'\oplus
\cF_{Y_2}''|_{\cD^{(1)}_{Y_2}\times X}$ whose image is the
diagonal $\cO$. Hence $\cE_{Y_2}|_{\cD^{(1)}_{Y_2}\times X}$ is a
family of extensions of a line bundle by rank 3 bundles.
 This extension splits
along $\tdel\cap Y_2$ so that we have an embedding of $\cO$ into
$\cE_{Y_2}|_{\tdel\cap Y_2\times X}$.

Note that $\tY$ is the blow-up of $Y_2$ along $\tdel\cap Y_2$ with
unstable points removed. Let $\cE_{\tY}$ be the pull-back of
$\cE_{Y_2}$ to $\tY\times X$ and $\cD^{(3)}_{\tY}$ be the
exceptional divisor while $\cD^{(1)}_{\tY}$ and $\cD^{(2)}_{\tY}$
denote the proper transforms of $\cD^{(1)}_{Y_2}$ and
$\cD^{(2)}_{Y_2}$ respectively.
 Let $\tcE$ be the kernel of
$$\cE_{\tY}\to \cE_{\tY}|_{\cD^{(3)}_{\tY}\times X}\to
\cE_{\tY}|_{\cD^{(3)}_{\tY}\times X}/\cO$$
 This is the desired family of semistable bundles of
rank 4. Verifying that this satisfies the conditions of
Proposition \ref{key} is a repetition of the computations in the
previous subsections and so we leave it to the reader.

\section{Blowing down Kirwan's desingularization}
In this section we show that the morphism
$$\rho:K\to S$$
constructed in section \ref{constr}, is in fact the result of two
contractions. In \cite{ogrady}, O'Grady worked out such
contractions for the moduli space of sheaves on a \emph{K3
surface}. We follow O'Grady's arguments to show that $K$ can be
contracted twice
\begin{equation}\label{eqf}\xymatrix{
f:&K\ar[r]^{f_{\sigma}}&K_{\sigma}\ar[r]^{f_{\epsilon}}&K_{\e}}\end{equation}
and these contractions are actually blow-downs. Then we show that
the map $\rho$ factors through $K_{\e}$, i.e.
\begin{equation}\label{facd}\xymatrix{
K\ar[rr]^{\rho}\ar[dr]_{f}&& S\\
&K_{\e}\ar[ur]_{\rho_{\e}}}\end{equation}
 By Zariski's main theorem, we
will conclude that $K_{\e}\cong S$.

\subsection{Contractions}\label{contractions}
Since the details are almost identical to section 3 of
\cite{ogrady}, we provide only the outline.

Let $\cA$ (resp. $\cB$) be the tautological rank 2 (resp. rank 3)
bundle over the Grassmannian $Gr(2,g)$ (resp. $Gr(3,g)$). Let
$W=sl(2)^{\vee}$ be the dual vector space of $sl(2)$. Fix $B\in
Gr(3,g)$. Then the variety of complete conics $\mathbf{CC} (B)$ is
the blow-up
$$\xymatrix{
\pp (S^2B)&\mathbf{CC}(B)\ar[l]_{\Phi_B}\ar[r]^{\Phi^{\vee}_B}&\pp
(S^2 B^{\vee})}$$ of both of the spaces of conics in $\pp B$ and
$\pp B^{\vee}$ along the locus of rank 1 conics.

\begin{proposition}\label{prop5.1}
\begin{enumerate}
\item $\tD^{(1)}$ is the variety of complete conics
$\mathbf{CC}(\cB)$ over $Gr(3,g)$. In other words, $\tD^{(1)}$ is
the blow-up of the projective bundle $\pp (S^2\cB)$ along the
locus of rank 1 conics.
\item There is an integer $l$ such that $$\tD^{(3)}\cong
\pp (S^2\cA)\times_{Gr(2,g)} \pp (\cc^g/\cA\oplus \cO(l)).$$ Hence
$\tD^{(3)}$ is a $\pp^2\times \pp^{g-2}$ bundle over $Gr(2,g)$.
\item The intersection $\tD^{(1)}\cap \tD^{(3)}$ is isomorphic to
the fibred product
$$\pp(S^2\cA)\times \pp (\cc^g/\cA)$$ over $Gr(2,g)$. As a
subvariety of $\tD^{(1)}$, $\tD^{(1)}\cap \tD^{(3)}$ is the
exceptional divisor of the blow-up $\mathbf{CC}(\cB)\to
{\pp}(S^2\cB^{\vee})$.
\item The intersection $\tD^{(1)}\cap\tD^{(2)}\cap \tD^{(3)}$ is isomorphic to
$$\pp(S^2\cA)_1\times \pp (\cc^g/\cA)$$ over $Gr(2,g)$ where $\pp(S^2\cA)_1$
denotes the locus of rank 1 quadratic forms.
\item The intersection $\tD^{(1)}\cap \tD^{(2)}$ is the
exceptional divisor of the blow-up $\mathbf{CC}(\cB)\to
{\pp}(S^2\cB)$.
\end{enumerate}\end{proposition}
\begin{proof}
The proofs are identical to (3.1.1), (3.5.1), and (3.5.4) in
\cite{ogrady}.\end{proof}

Next, we consider some rational curves to be contracted. Define
the following classes in $N_1(\tD^{(1)})$ (the group of numerical
equivalence classes of 1-cycles)
$$\s:=\text{ the class of lines in the fiber of }\Phi_B^{\vee}$$
$$\e:=\text{ the class of lines in the fiber of }\Phi_B$$
$$\gamma:=\text{ the class of }\{\Phi_{B_t}^{-1}(q_t)\}_{t\in \Lambda}$$
where $\{B_t\}$ is a line $\Lambda$ of 3-dimensional subspaces in
$Gr(3,g)$ containing a fixed 2-dimensional space $A$ with $q\in
S^2A$ and $q_t$ is the induced quadratic form on $B_t$.

To show that these form a basis of $N_1(\tD^{(1)})$ we consider
the following  diagram
$$\xymatrix{ \tD^{(1)}\ar[r]^{\theta}&\pp(S^2\cB)\ar[d]^{\phi}\\
&Gr(3,g)}$$ where $\theta$ is the blow-up. Let
$h=c_1(\cB^{\vee})$, $x=c_1(\cO_{\pp (S^2\cB)}(1))$ and $e$ be the
exceptional divisor of $\theta$. Then obviously $h,x,e$ form a
basis of $N^1(\tD^{(1)})$ which is dual to $N_1(\tD^{(1)})$. By
elementary computation, the intersection pairing is given by the
table
$$\begin{matrix}
& h& x&e\\
\e&0&0&-1\\
\s&0&1&2\\
\gamma&1&0&0\end{matrix}$$ Hence, $\s,\e,\gamma$ form a basis of
$N_1(\tD^{(1)})$.

\begin{lemma}\begin{enumerate}
\item $[\tD^{(1)}]|_{\mathbf{CC}(B)}=-2x+e|_{\mathbf{CC}(B)}$ for $B\in Gr(3,g)$.
\item $[\tD^{(2)}]|_{\tD^{(1)}}=e$
\item $[\tD^{(3)}]|_{\tD^{(1)}}=3x-2h-2e$
\item $\Theta_{\tD^{(1)}}=-(g-4)h-6x+2e$ where
$\Theta_{\tD^{(1)}}$ denotes the canonical divisor of $\tD^{(1)}$.
\end{enumerate}\end{lemma}
The proofs are identical to those of (3.2.3) - (3.2.5), (3.4.3)
with obvious modifications.

Let $\hat{\s}=\imath_*\s$, $\he=\imath_*\e$ and $\hg=\imath_*\g$
where $\imath$ is the inclusion of $\tD^{(1)}$ into $K$. By the
above lemma, $x, h, e$ are in the image of $N^1(K)$ by
restriction. Hence, $N^1(K)\to N^1(\tD^{(1)})$ is surjective and
dually $\imath_*$ is injective. Consequently,
 $\hs, \he, \hg$ are linearly independent.

At this point, we can compute the discrepancy
$\omega_{K}-\pi^*\omega_{M_0}$ of the canonical divisors
$\omega_{K}$ and $\omega_{M_0}$.

\begin{proposition}\label{discr}
$$\omega_{K}-\pi^*\omega_{M_0}=(3g-1)\tD^{(1)}+(g-2) \tD^{(2)}+(2g-2)
\tD^{(3)}$$
\end{proposition}
\begin{proof} Obvious adaptation of the proof of (3.4.1) in \cite{ogrady}.
\end{proof}

\begin{corollary}\label{cor5.4} For $g\ge 3$, $M_0$ has terminal singularities and the
plurigenera are all trivial.\end{corollary}
\begin{proof} It is well-known that $\omega_{M_0}$ is anti-ample.
Since the singularities are terminal,
$\pi_*\omega_K=\omega_{M_0}$. It follows from spectral sequence
and Kodaira's vanishing theorem that $H^0(K, \omega_K^{\otimes
m})\cong H^0(M_0, \omega_{M_0}^{\otimes m})=0$ for
$m>0$.\end{proof}

Finally we can show that $K$ can be blown-down twice.
\begin{proposition}\begin{enumerate}
\item $\hs, \he$ are $\omega_K$-negative extremal rays. For $g>3$,
$\hg$ is also $\omega_K$-negative extremal.
\item The contraction $K_{\s}$ of the ray $\rr^+\hs$ is a smooth
projective desingularization of $M_0$. In fact, this is the
contraction of the $\pp(S^2\cA)$-direction of $\tD^{(3)}$. Since
the normal bundle is $\cO(-1)$ up to tensoring a line bundle on
$\pp (\cc^g/\cA\oplus \cO(l))$, the contraction is a blow-down
map.
\item The image of $\he$ in $N_1(K_{\s})$ is
$\omega_{K_{\s}}$-negative extremal ray and its contraction
$K_{\e}$ is a smooth projective desingularization of $M_0$. This
is the contraction of the fiber direction of
$\pp(S^2\cB^{\vee})\to Gr(3,g)$ and is also a blow-down map.
\end{enumerate}\end{proposition}
The proofs are same as those of (3.0.2)-(3.0.4) in \cite{ogrady}.

\subsection{Factorization of $\rho$}
Now we can show the following
\begin{theorem}\label{th5.6} $\rho$ factors through $K_{\e}$ and $K_{\e}\cong
S$.\end{theorem}
\begin{proof}
Let us consider the first contraction $f_{\s}:K\to K_{\s}$. We
claim that there is a continuous map $\rho_{\s}: K_{\s}\to S$ such
that $\rho_{\s}\circ f_{\s}=\rho$. (See the diagram \eqref{facd}.)
By Riemann's extension theorem \cite{Kran}, it suffices to show
that $\rho$ is constant on the fibers of $f_{\s}$. From
Proposition \ref{prop5.1}, we know $f_{\s}$ is the result of
contracting the fibers $\pp^2$  of
$$\tD^{(3)}=\pp (S^2\cA)\times \pp (\cc^g/\cA\oplus \cO(l))\to
\pp (\cc^g/\cA\oplus \cO(l))$$ which amounts to forgetting the
choice of $b, c$ in the 2-dimensional subspace of $H^1(\cO)$
spanned by $b,c$. We need only to check that the isomorphism
classes of the vector bundles given by \eqref{eq4.9} and
\eqref{eq4.10} depend \emph{not} on the particular choice of $b,c$
\emph{but} only on the points in $\pp^{g-2}$-bundle $\pp
(\cc^g/\cA\oplus \cO(l))\to \pp (\cc^g/\cA\oplus \cO(l))$ over
$Gr(2,g)$.

From \cite{BaSe} Proposition 5, the isomorphism classes of bundles
given by \eqref{eq4.9} are parametrized by a vector bundle of rank
$g-2$  over $Gr(2,g)$. In particular, the isomorphism classes are
independent of the choice of $b,c$.  Hence the bundles given by
\eqref{eq4.9} are constant along the $\pp(S^2\cA)$-direction.  On
the other hand, it is elementary to show that a similar statement
holds for the bundles given by \eqref{eq4.10}. Therefore, there
exists a morphism $\rho_{\s}:K_{\s}\to S$ such that
$\rho_{\s}\circ f_{\s}=\rho$.

Next we show that $\rho_{\s}$ factors through $K_{\e}$.  The
morphism $f_{\e}:K_{\s}\to K_{\e}$ is the contraction of the
fibers $\pp^5$ of
$$ \pp (S^2\cB)\to Gr(3,g)$$
 and general points
of a fiber give rise to a rank 4 bundle whose transition matrices
are of the form \eqref{eq4.5}. It is elementary to show that the
isomorphism classes of the bundles given by \eqref{eq4.5} depend
only on the 3-dimensional subspace spanned by $a,b,c$. Hence
$\rho_{\s}$ is constant along the fibers of $f_{\e}$. By Riemann's
extension theorem again, we get a morphism $\rho_{\e}:K_{\e}\to S$
such that $\rho_{\e}\circ f=\rho$.

From \cite{Ba, BaSe}, $\rho(\tD^{(2)}-\tD^{(1)}\cup \tD^{(3)})$ is
a smooth divisor of $S-\rho(\tD^{(1)}\cup \tD^{(3)})$ that lies
over $\mathfrak{K}-\zz_2^{2g}$. Hence, we have a morphism from
$S-\rho(\tD^{(1)}\cup \tD^{(3)})$ to the blow-up of
$M_0-\zz_2^{2g}$ along $\mathfrak{K}-\zz_2^{2g}$ which is
isomorphic to $K-\tD^{(1)}\cup \tD^{(3)}=K_{\e}-f(\tD^{(1)}\cup
\tD^{(3)})$ by construction. Hence, $\rho_{\e}$ is an isomorphism
in codimension one. Since $K_{\e}$ and $S$ are both smooth,
Zariski's main theorem says $K_{\e}$ is isomorphic to $S$.
\end{proof}

\begin{conjecture}
The intermediate variety $K_{\s}$ is the Narasimhan-Ramanan
desingularization.\end{conjecture}
 We hope to get back to this
conjecture in the future.

\section{Cohomological consequences}

\subsection{Cohomology of Seshadri's desingularization}
In \cite{Ba, BaSe}, Balaji and Seshadri show the Betti numbers of
Seshadri's desingularization $S$ can be computed, up to degree
$\le 2g-4$. Thanks to the explicit description of $S$ as the
blow-down of $K$, we can compute the Betti numbers in all degrees.

For a variety $T$, let $$P(T)=\sum_{k=0}^{\infty} t^k\dim H^k(T)$$
be the Poincar\'e series of $T$. In \cite{k2}, Kirwan described an
algorithm for the Poincar\'e series of a partial desingularization
of a good quotient of a smooth projective variety and in \cite{k5}
the algorithm was applied to the moduli space without fixing the
determinant. For $P(M_2)$ we use Kirwan's algorithm in \cite{k2}.

By \cite{AB} \S11 and \cite{KirL}, it is well-known that the
equivariant Poincar\'e series $P^G(\frak{R}^{ss})=\sum_{k\ge 0}t^k
\dim H^k_G(\frak{R}^{ss})$ is
$$P^G(\frak{R}^{ss})=\frac{(1+t^3)^{2g}-t^{2g+2}(1+t)^{2g}}{(1-t^2)(1-t^4)}$$
up to degrees as high as we want. In order to get
$\frak{R}^{ss}_1$ we blow up $\frak{R}^{ss}$ along
$GZ^{ss}_{SL(2)}$ and delete the unstable strata. So we get
$$P^G(\frak{R}_1^{ss})=P^G(\frak{R}^{ss})+2^{2g}\big(\frac{t^2+t^4+\cdots+t^{6g-2}}{1-t^4}
-\frac{t^{4g-2}(1+t^2+\cdots+t^{2g-2})}{1-t^2}\big).$$ Now
$\frak{R}_2^{ss}$ is obtained by blowing up $\frak{R}_1^{ss}$
along $G\tilde{Z}^{ss}_{\cc^*}$ and deleting the unstable strata.
Thus we have
\begin{equation}\begin{array}{ll}
P^G(\frak{R}_2^{ss})=P^G(\frak{R}_1^{ss})
&+(t^2+t^4+\cdots+t^{4g-6})\big(\frac12
\frac{(1+t)^{2g}}{1-t^2}+\frac12
\frac{(1-t)^{2g}}{1+t^2} +2^{2g}\frac{t^2+\cdots+t^{2g-2}}{1-t^4}\big)\\
&-\frac{t^{2g-2}(1+t^2+\cdots+t^{2g-4})}{1-t^2}\big((1+t)^{2g}+2^{2g}(t^2+t^4+\cdots+t^{2g-2})\big).
\end{array}
\end{equation}
Because the stabilizers of the $G$ action on $\frak{R}^{ss}_2$ are
all finite, we have $$H^*_G(\frak{R}_2^{ss})\cong
H^*(\frak{R}^{ss}_2/G)=H^*(M_2)$$ and hence we deduce that
\begin{equation}\label{N2Pr}
\begin{array}{ll}
P(M_2)&=\frac{(1+t^3)^{2g}-t^{2g+2}(1+t)^{2g}}{(1-t^2)(1-t^4)}\\
&+2^{2g}\big(\frac{t^2+t^4+\cdots+t^{6g-2}}{1-t^4}
-\frac{t^{4g-2}(1+t^2+\cdots+t^{2g-2})}{1-t^2}\big)\\
&+(t^2+t^4+\cdots+t^{4g-6})\big(\frac12
\frac{(1+t)^{2g}}{1-t^2}+\frac12
\frac{(1-t)^{2g}}{1+t^2} +2^{2g}\frac{t^2+\cdots+t^{2g-2}}{1-t^4}\big)\\
&-\frac{t^{2g-2}(1+t^2+\cdots+t^{2g-4})}{1-t^2}\big((1+t)^{2g}+2^{2g}(t^2+t^4+\cdots+t^{2g-2})\big).
\end{array}
\end{equation}
Kirwan's desingularization is the blow-up of $M_2$ along
$\tilde{\Delta}\git SL(2)$ which is isomorphic to the $2^{2g}$
copies of  $\pp(S^2\cA)$ over $Gr(2,g)$. Hence,
$$P(K)=P(M_2)+2^{2g}(1+t^2+t^4)P(Gr(2,g))(t^2+t^4+\cdots + t^{2g-4})$$
by \cite{GH} p. 605.\footnote{The formula in \cite{GH} is stated
for smooth manifolds. But the same Mayer-Vietoris argument gives
us the same formula in our case (of orbifold $M_2$ blown up along
a smooth subvariety). The only thing to be checked is that the
pull-back homomorphism $H^*(M_2)\to H^*(K)$ is injective but this
clearly holds by the decomposition theorem of Beilinson,
Bernstein, Deligne and Gabber.}

On the other hand, $K$ is the blow-up of $K_{\s}$ along a
$\pp^{g-2}$-bundle over $Gr(2,g)$. Hence,
$$\begin{array}{ll}
P(K_{\s})&=P(K)-2^{2g}(1+t^2+\cdots+t^{2g-4})P(Gr(2,g))(t^2+t^4)\\
&=P(M_2)+2^{2g}P(Gr(2,g))\frac{t^6-t^{2g-2}}{1-t^2}.\end{array}$$
Similarly, $K_{\s}$ is the blow-up of $K_{\e}$ along a $Gr(3,g)$
and thus
$$\begin{array}{ll}
P(K_{\e})&=P(K_{\s})-2^{2g}P(Gr(3,g))(t^2+\cdots+t^{10})\\
&=P(M_2)+2^{2g}P(Gr(2,g))\frac{t^6-t^{2g-2}}{1-t^2}-2^{2g}P(Gr(3,g))(t^2+\cdots+t^{10}).
\end{array}$$
Since $K_{\e}$ is isomorphic to Seshadri's desingularization, we
get
$$\begin{array}{ll}
P(S)&=\frac{(1+t^3)^{2g}-t^{2g+2}(1+t)^{2g}}{(1-t^2)(1-t^4)}\\
&+2^{2g}\big(\frac{t^2+t^4+\cdots+t^{6g-2}}{1-t^4}
-\frac{t^{4g-2}(1+t^2+\cdots+t^{2g-2})}{1-t^2}\big)\\
&+(t^2+t^4+\cdots+t^{4g-6})\big(\frac12
\frac{(1+t)^{2g}}{1-t^2}+\frac12
\frac{(1-t)^{2g}}{1+t^2} +2^{2g}\frac{t^2+\cdots+t^{2g-2}}{1-t^4}\big)\\
&-\frac{t^{2g-2}(1+t^2+\cdots+t^{2g-4})}{1-t^2}\big((1+t)^{2g}+2^{2g}(t^2+t^4+\cdots+t^{2g-2})\big)\\
&+2^{2g}P(Gr(2,g))\frac{t^6-t^{2g-2}}{1-t^2}-2^{2g}P(Gr(3,g))(t^2+\cdots+t^{10}).
\end{array}$$
By Schubert calculus \cite{GH}, we have
$$P(Gr(2,g))=\frac{(1-t^{2g})(1-t^{2g-2})}{(1-t^2)(1-t^4)}$$
$$P(Gr(3,g))=\frac{(1-t^{2g})(1-t^{2g-2})(1-t^{2g-4})}{(1-t^2)(1-t^4)(1-t^6)}$$
and hence we obtained a closed formula for the Poincar\'e
polynomial of $S$.

In \cite{BaSe}, an algorithm for the Betti numbers only up to
degree $2g-4$ is provided. It is an elementary exercise to check
that in this range, their answer is identical to ours.

\subsection{The stringy E-function}
The stringy E-function is an invariant of singular varieties
introduced by Batyrev, Denef and Loeser, based on the suggestions
by Kontsevich. In \cite{kiem}, the stringy E-function of $M_0$ was
computed for $g=3$ by using the observation that the singularities
are hypersurface singularities in this case.\footnote{There is a
small error in \cite{kiem} page 1852. In line -3, $\alpha_1$
should be replaced by $\alpha_7^2$ and thus in line -1, the
discrepancy divisor is $8D_1+D_2+4D_3$ (cf. Proposition
\ref{discr}). The computation in \cite{kiem} \S7 should be
accordingly modified. The correct formula for any $g\ge 3$ is
proved in this paper (Theorem \ref{th6.1}).} In this subsection,
we compute the stringy E-function of $M_0$ for arbitrary genus.
For the definition and some basic facts on the stringy
E-functions, see the introduction of \cite{kiem}.

Since the discrepancy divisor is given by Proposition \ref{discr},
our goal is to compute
$$\begin{array}{ll}
E_{st}(M_0)&=E(M^s_0)+E(\tD^{(1)}_0)\frac{uv-1}{(uv)^{3g}-1} +
E(\tD^{(2)}_0)\frac{uv-1}{(uv)^{g-1}-1}+
E(\tD^{(3)}_0)\frac{uv-1}{(uv)^{2g-1}-1}\\
&+E(\tD^{(1,2)}_0)\frac{uv-1}{(uv)^{3g}-1}\frac{uv-1}{(uv)^{g-1}-1}+
E(\tD^{(2,3)}_0)\frac{uv-1}{(uv)^{g-1}-1}\frac{uv-1}{(uv)^{2g-1}-1}\\
&+E(\tD^{(1,3)}_0)\frac{uv-1}{(uv)^{3g}-1}\frac{uv-1}{(uv)^{2g-1}-1}+
E(\tD^{(1,2,3)}_0)\frac{uv-1}{(uv)^{3g}-1}\frac{uv-1}{(uv)^{g-1}-1}\frac{uv-1}{(uv)^{2g-1}-1}
\end{array}$$
where $\tD^{(I)}_0=\cap_{i\in I}\tD^{(i)}-\cup_{j\notin
I}\tD^{(j)}$ for $I\subset \{1,2,3\}$ and $E$ denotes the
Hodge-Deligne polynomal.

The E-function of the smooth part is from \cite{kiem} \S 4,
$$\begin{array}{ll}
E(M^s_0)&=E(M_2)-E(D_2^{(1)})-E(D_2^{(2)}-D_2^{(1)})\\
&=\frac{(1-u^2v)^g(1-uv^2)^g-(uv)^{g+1}(1-u)^g(1-v)^g}{(1-uv)(1-(uv)^2)}\\
&-\frac12(\frac{(1-u)^g(1-v)^g}{1-uv}+\frac{(1+u)^g(1+v)^g}{1+uv}).
\end{array}$$

By Proposition \ref{prop5.1},
$\tD^{(1)}_0=\tD^{(1)}-(\tD^{(2)}\cup \tD^{(3)})$ is the union of
$2^{2g}$ copies of $\pp^5-\pp^2\times_{\zz_2}\pp^2$-bundle over
$Gr(3,g)$ and thus
$$E(\tD^{(1)}_0)\frac{uv-1}{(uv)^{3g}-1}=2^{2g}((uv)^5-(uv)^2)E(Gr(3,g))
\frac{uv-1}{(uv)^{3g}-1}.$$

Since $\tD^{(2)}_0$ is the quotient of a $\pp^{g-2}\times
\pp^{g-2}$-bundle over $Jac_0-\zz_2^{2g}$ by the action of
$\zz_2$, the E-function of $\tD^{(2)}_0$ is
$$\begin{array}{ll}
&E(\tD^{(2)}_0)\frac{uv-1}{(uv)^{g-1}-1} \\
&=\big(\frac12 (1-u)^g(1-v)^g+\frac12 (1+u)^g(1+v)^g
-2^{2g}\big)E(\pp^{g-2}\times \pp^{g-2})^+\frac{uv-1}{(uv)^{g-1}-1}\\
&+\big(\frac12 (1-u)^g(1-v)^g-\frac12 (1+u)^g(1+v)^g
\big)E(\pp^{g-2}\times \pp^{g-2})^-
\frac{uv-1}{(uv)^{g-1}-1}\end{array}$$ where
$$E(\pp^{g-2}\times \pp^{g-2})^+=\frac{((uv)^g-1)((uv)^{g-1}-1)}{(uv-1)((uv)^2-1)}$$
is the E-polynomial of the $\zz_2$-invariant part of
$H^*(\pp^{g-2}\times \pp^{g-2})$ and
$$E(\pp^{g-2}\times \pp^{g-2})^-=uv\frac{((uv)^{g-1}-1)((uv)^{g-2}-1)}{(uv-1)((uv)^2-1)}$$
is the E-polynomial of the anti-invariant part.

By Proposition \ref{prop5.1}, $\tD^{(3)}_0$ is the union of
$2^{2g}$ copies of a $(\pp^2\times \pp^{g-2}-\pp^2\times
\pp^{g-3}\cup \pp^1\times \pp^{g-2})$-bundle over $Gr(2,g)$ and
thus
$$E(\tD^{(3)}_0)\frac{uv-1}{(uv)^{2g-1}-1}=2^{2g}(uv)^gE(Gr(2,g))\frac{uv-1}{(uv)^{2g-1}-1}.$$

Notice that $\tD^{(1,2)}_0$ is the disjoint union of $2^{2g}$
copies of a $(\pp^2-\pp^1)\times \pp^2$-bundle over $Gr(3,g)$ and
thus
$$E(\tD^{(1,2)}_0)\frac{uv-1}{(uv)^{3g}-1}\frac{uv-1}{(uv)^{g-1}-1}=
2^{2g}((uv)^2+(uv)^3+(uv)^4)E(Gr(3,g))\frac{uv-1}{(uv)^{3g}-1}\frac{uv-1}{(uv)^{g-1}-1}.$$

Also, $\tD^{(1,3)}_0$ is a $(\pp^2-\pp^1)\times \pp^{g-3}$-bundle
over $Gr(2,g)$ and thus
$$E(\tD^{(1,3)}_0)\frac{uv-1}{(uv)^{3g}-1}\frac{uv-1}{(uv)^{2g-1}-1}=
2^{2g}(uv)^2\frac{(uv)^{g-2}-1}{uv-1}E(Gr(2,g))\frac{uv-1}{(uv)^{3g}-1}\frac{uv-1}{(uv)^{2g-1}-1}.$$

Finally, a component of $\tD^{(2,3)}_0$ is a $\pp^1\times
(\pp^{g-2}-\pp^{g-3})$-bundle over $Gr(2,g)$ and a component of
$\tD^{(1,2,3)}_0$ is a $\pp^1\times \pp^{g-3}$-bundle over
$Gr(2,g)$. Therefore,
$$E(\tD^{(2,3)}_0)\frac{uv-1}{(uv)^{g-1}-1}\frac{uv-1}{(uv)^{2g-1}-1}=
2^{2g}(1+uv)(uv)^{g-2}E(Gr(2,g))\frac{uv-1}{(uv)^{g-1}-1}\frac{uv-1}{(uv)^{2g-1}-1}$$
and
$$\begin{array}{ll}E(\tD^{(1,2,3)}_0)&\frac{uv-1}{(uv)^{3g}-1}
\frac{uv-1}{(uv)^{g-1}-1}\frac{uv-1}{(uv)^{2g-1}-1}\\&=
2^{2g}(1+uv)\frac{(uv)^{g-2}-1}{uv-1}E(Gr(2,g))\frac{uv-1}{(uv)^{3g}-1}
\frac{uv-1}{(uv)^{g-1}-1}\frac{uv-1}{(uv)^{2g-1}-1}.\end{array}$$
Recall that
$$E(Gr(2,g))=\frac{((uv)^g-1)((uv)^{g-1}-1)}{(uv-1)((uv)^2-1)}$$
$$E(Gr(3,g))=\frac{((uv)^g-1)((uv)^{g-1}-1)((uv)^{g-2}-1)}{(uv-1)((uv)^2-1)((uv)^3-1)}.$$

Putting together all the pieces above, we get
\begin{theorem}\label{th6.1}
$$\begin{array}{ll}E_{st}(M_0)=&
\frac{(1-u^2v)^g(1-uv^2)^g-(uv)^{g+1}(1-u)^g(1-v)^g}{(1-uv)(1-(uv)^2)}\\
&- \frac{(uv)^{g-1}}{2}\big(
\frac{(1-u)^g(1-v)^g}{1-uv}-\frac{(1+u)^g(1+v)^g}{1+uv} \big).
\end{array}$$\end{theorem}

\begin{remark} It is well-known that the middle perversity
intersection cohomology of $M_0$ is
equipped with a Hodge structure and hence it makes sense to think
about the E-polynomial of the intersection cohomology. The
computation of the Poincar\'e polynomial of $IH^*(M_0)$ in
\cite{k5} can be easily refined as in \cite{EK} to give the
E-polynomial of $IH^*(M_0)$
$$\begin{array}{ll}
IE(M_0)&=\frac{(1-u^2v)^g(1-uv^2)^g-(uv)^{g+1}(1-u)^g(1-v)^g}{(1-uv)(1-(uv)^2)}\\
&- \frac{(uv)^{g-1}}{2}\big(
\frac{(1-u)^g(1-v)^g}{1-uv}+(-1)^{g-1}\frac{(1+u)^g(1+v)^g}{1+uv}
\big).\end{array}
$$ See also \cite{kiem2}.
Quite surprisingly, when $g$ is even, $E_{st}(M_0)$ is identical
to the E-polynomial of the middle perversity intersection
cohomology of $M_0$. This indicates that there may be an unknown
relation between the stringy E-function and the intersection
cohomology. When $g$ is odd, $E_{st}(M_0)$ is not a polynomial.
\end{remark}

\begin{corollary}
The stringy Euler number of $M_0$ is $$e_{st}(M_0):=\lim_{u,v\to
1}E_{st}(M_0)=4^{g-1}.$$\end{corollary}

Let $e_g$ be the stringy Euler number of the moduli space $M_0$
for a genus $g$ curve. When $g=2$, $M_0\cong \pp^3$ and so
$e_2=4$. Therefore the equality
$$\sum_{g}e_gq^g=\frac14\frac1{1-4q}$$
holds for degree $\ge 2$. The coefficient $\frac14$ might be
related to the ``mysterious" coefficient $\frac14$ for the
S-duality conjecture test in the K3 case in \cite{VW}.


\end{document}